\documentstyle[amssymb,amstex,12pt,a4wide]{article}
\newtheorem{theorem}{Theorem}

\newtheorem{proposition}{Proposition}
\newtheorem{corollary}{Corollary}

\newcommand{\vol}{\operatorname{Vol}}
\newcommand{\pl}{\operatorname{pl}}
\newcommand{\Area}{\operatorname{Area}}
\newcommand{\entry}{\operatorname{Entry}}
\newcommand{\innt}{\operatorname{Int}}
\begin{document}
\setlength{\baselineskip}{1.7\baselineskip}
\title{A polyhedral Markov field - pushing the Arak-Surgailis
       construction into three dimensions}
\author{Tomasz Schreiber\footnote{Research supported by the Foundation for
 Polish Science (FNP), by the Polish Minister of
 Scientific Research and Information Technology grant 1 P03A 018 28
 (2005-2007) and from
 the EC 6th Framework Programme Priority 2 Information Society Technology
 Network of Excellence MUSCLE  (Multimedia Understanding through Semantics,
 Computation and Learning; FP6-507752). A part of this research was completed
 while staying at the Centrum voor Wiskunde en Informatica (CWI), Amsterdam,
 The Netherlands}\\
        Faculty of Mathematics \& Computer Science,\\
        Nicolaus Copernicus University,\\
        ul. Chopina 12 $\slash$ 18,\\
        Toru\'n, Poland,\\
        e-mail: tomeks at mat.uni.torun.pl}
\date{}
\maketitle

\paragraph{Abstract:} {\it The purpose of the paper is to construct
 a polyhedral Markov field in ${\Bbb R}^3$ in analogy with the
 planar construction of the original Arak (1982) polygonal
 Markov field. We provide a dynamic construction of the process
 in terms of evolution of two-dimensional multi-edge
 systems tracing polyhedral boundaries of the field in
 three-dimensional time-space. We also give a general
 algorithm for simulating Gibbsian modifications
 of the constructed polyhedral field. }

\section{Introduction}

 The notion of a consistent polygonal Markov field has first appeared in the
 seminal paper by Arak (1982) who constructed an isometry-invariant process
 with polygonal realisations in the plane, enjoying a two-dimensional
 germ Markov property [the conditional behaviour of the field in an
 open bounded domain with piecewise smooth boundary depends on the
 exterior configuration only through arbitrarily close neighbourhoods
 of the boundary] and with one-dimensional sections coinciding in
 distribution with homogeneous Poisson point processes. An attractive
 feature of this process is that it admits several alternative equivalent
 representations. Two of them fall into the general Gibbsian framework
 and use, respectively, lines and points as the basic building blocks,
 see Arak \& Surgailis (1989), Arak \& Surgailis (1991) and Arak,
 Clifford \& Surgailis (1993). The third representation of the Arak-Surgailis
 field is available in terms of equilibrium evolution of one-dimensional
 particle systems, tracing the polygonal realisations of the process in
 two-dimensional time-space, see the above papers. A dynamic representation in this
 spirit is in fact  shared by a much richer class of planar polygonal Markov
 processes, as constructed and discussed in Arak \& Surgailis (1989,1991)
 and Arak, Clifford \& Surgailis (1993).

 The purpose of the present article is to construct an analogue of
 the Arak process in three-dimensional polyhedral setting. We
 provide a plane-based Gibbsian representation as well as
 a dynamic construction in terms of evolving multi-edge systems
 tracing polyhedral boundaries of the field in three-dimensional
 time-space ${\Bbb R} \times {\Bbb R}^2.$ The resulting field
 is Markovian, as directly follows by the Gibbsian representation.
 The construction is isometry invariant. We give also an explicit
 formula for the corresponding partition function. However, we are
 not able to explicitly characterise the equilibrium initial condition
 for the constructed dynamics for multi-edge systems and hence to
 explicitly construct a consistent version of our polyhedral Markov
 field.

 The results discussed above are presented Sections \ref{KOPROC}
 and \ref{WLAA}. They are complemented with Section \ref{BSDYN} where
 we provide a Metropolis type simulation algorithm for a general
 class of Gibbsian modifications of the polyhedral Markov field
 constructed in this paper. Apart from its theoretical interest
 the purpose of developing this algorithm is related to the fact
 that we anticipate possible applications of polyhedral Markov
 fields for volumetric image segmentation, in the spirit of
 our papers Kluszczy\'nski et. al. (2004,2005), where similar
 algorithms in two-dimensional setting are used for planar image
 segmentation, see also Clifford \& Middleton (1989) and
 Clifford \& Nicholls (1994).

\section{Construction of the process}\label{KOPROC}

\subsection{Preliminaries}
 In what follows we abuse the language granting the name {\it polygonal face}
 to any open connected subset of a plane in ${\Bbb R}^3,$ with polygonal
 boundary. In particular, a {\it face} may well contain holes but it
 cannot split into several disconnected parts. For an open bounded
 and convex polyhedral domain $D \subseteq {\Bbb R}^3$ we define the
 family $\Gamma_D$ of admissible polyhedral configurations in $D$ by
 taking all collections $\gamma$ of planar polygonal faces contained
 in $D$ such that
 \begin{description}
  \item{\bf (P1)} the faces of $\gamma$ do not intersect, but they may share edges,
  \item{\bf (P2)} each interior edge of $\gamma$ (contained in $D$) is shared by
                  exactly two faces,
  \item{\bf (P3)} each boundary edge of $\gamma$ (contained in $\partial D$)
                  belongs to exactly one face,
  \item{\bf (P4)} each internal vertex of $\gamma$ (contained in $D$) is
                  shared by exactly three faces and has exactly three
                  outgoing edges,
  \item{\bf (P5)} each boundary vertex of $\gamma$ (contained in $\partial D$)
                  is shared by exactly two faces and has exactly one outgoing
                  internal edge and exactly two outgoing boundary edges,
  \item{\bf (P6)} no two faces of $\gamma$ are coplanar.
 \end{description}
  In other words, an admissible polyhedral configuration $\gamma$ corresponds
  to a family of disjoint closed polyhedral surfaces in $D,$ possibly nested
  and chopped off by the boundary $\partial D.$ We write $E_D(\gamma)$ for the
  collection of edges of $\gamma$ contained in $D$ (internal edges),
  $E_{\partial D}(\gamma)$ for the collection of edges of $\gamma$
  contained in $\partial D$ (boundary edges), and $F_D(\gamma)$ for the
  collection of faces of $\gamma$ contained in $D$ (only internal faces
  are considered). For a finite collection
  $\{ \varpi_i \}_{i=1}^n$ of planes intersecting $D$ we consider
  the family $\Gamma_D(\{ \varpi_i \}_{i=1}^{n})$ of all admissible
  polyhedral configurations $\gamma$ in $D$ with the additional properties
  that $\gamma \subseteq \bigcup_{i=1}^n \varpi_i$ and that $\gamma \cap
   \varpi_i$ is a single polygonal face of non-zero area, possibly with
  some zero-measure polygonal curves added, for each $\varpi_i, \;
  i=1,2,...,n.$ Let $\mu$ be the usual isometry-invariant Haar-Lebesgue
  measure on the space ${\cal P}$ of all two-dimensional planes in
  ${\Bbb R}^3.$ One possible construction of $\mu$ goes by identifying
  a plane $\varpi$ with the pair $(u,\rho) \in {\Bbb S}_2 \times {\Bbb R}_+,$
  with ${\Bbb S}_2$ standing for the unit sphere in ${\Bbb R}^3$ and
  where $u$ is the direction of the vector orthogonal to $\varpi$ and
  joining it to the origin, while $\rho$ is the distance between
  $\varpi$ and the origin. In these terms $\mu$ arises by endowing
  ${\Bbb S}_2 \times {\Bbb R}_+$ with the product of the usual
  surface measure on the sphere and the Lebesgue measure on
  ${\Bbb R}_+.$ We write $\Pi$ for the Poisson plane process with
  intensity measure $\mu$ and for each bounded domain $D$ we let $\Pi_D$
  stand for the restriction of $\Pi$ to the family ${\cal P}_D$ of planes
  intersecting $D.$ The following properties of the plane process
  $\Pi$ will be of use in the sequel. Here and throughout, to avoid
  possible confusion we use $d_2 x$ and $d_3 x$ rather than $dx$ to
  denote integration respectively w.r.t. the 2- and 3-dimensional
  Lebesgue measure, thus explicitly indicating the dimensionality
  of the integration variable.
  \begin{proposition}\label{WLASNOSCI}
   With the above notation we have:
   \begin{description}
    \item{\bf (I1)} The intersection of $\Pi$ with a given straight
          line $l$ in ${\Bbb R}^3$ is a homogeneous Poisson point
          process of intensity $\pi,$
    \item{\bf (I2)} The intersection of $\Pi$ with a given plane
          $\varpi$ is a Poisson process of lines $l = l(\phi,r)$
          in $\varpi$ with intensity measure $\frac{\pi}{2} d\phi dr,$
          where $r$ stands for the distance between $l$ and a certain
          fixed point $\hat{0}$ in $\varpi,$ while $\phi$ is the
          angle between some fixed line in $\varpi$ and the vector
          in $\varpi$ orthogonal to $l$ and joining it to $\hat{0},$
    \item{\bf (I3)} For a given plane $\varpi$ and $x \in \varpi$
          the probability that two planes of $\Pi$ meet at $[x,x+d_2x]
          \subset \varpi$ is $\frac{\pi^3}{4} d_2x.$
    \item{\bf (I4)} For $x \in {\Bbb R}^3$ the probability that three
          planes of $\Pi$ meet at $[x,x+d_3 x] \subset {\Bbb R}^3$ is
          $\frac{\pi^4}{6} d_3 x.$
  \end{description}
 \end{proposition}
  Indeed, the assertions {\bf (I1)} and {\bf (I2)} follow immediately by (3.29T)
  in  Miles (1971) because our process $\Pi$ coincides with ${\cal B}(2\pi,2,3)$
  there. Further, {\bf (I4)} is a direct consequence of (6-1-9') in Matheron (1975)
  with $a=\pi \slash 2$ there. Finally, {\bf (I3)} follows from {\bf (I2)}
  and (6-1-9') ibidem with $a = \pi \slash 2$ there.

 \subsection{Gibbsian representation}
  For a given admissible polyhedral configuration $\gamma$ we
  consider its {\it energy}
  \begin{equation}\label{ENERGIA}
    \Phi_D(\gamma) = 
       \frac{1}{2} \sum_{e \in E_D(\gamma)} (2\pi-|\angle(e)|) \ell(e)
    +  \frac{\pi^3}{4} \sum_{f \in F_D(\gamma)} \Area(f) + \frac{\pi^4}{6} \vol(D),
  \end{equation}
  where $\angle(e)$ is the solid convex angle between the planes determined
  by the faces meeting at an internal edge $e,$ $\ell(e)$ is the length
  of the edge $e,$ $\Area(f)$ stands for the area of the polygonal face
  $f$ while $\vol(D)$ denotes the volume of the domain $D.$
  We define our three-dimensional polyhedral process ${\cal A}_D$ in
  $D$ as the Gibbsian modification of the process induced on $\Gamma_D$
  by $\Pi_D$ with the Hamiltonian $\gamma \mapsto \Phi_D(\gamma).$
  To be more specific, we put
  \begin{equation}\label{DEFI1}
   {\Bbb P}\left({\cal A}_D \in {\cal F}\right)
   = \frac{{\Bbb E}\sum_{\gamma \in \Gamma_D(\Pi_D) \cap {\cal F}}
           \exp(-\Phi_D(\gamma))}
          {{\Bbb E}\sum_{\gamma \in \Gamma_D(\Pi_D)} \exp(-\Phi_D(\gamma))}
  \end{equation}
  for all ${\cal F} \subseteq \Gamma_D$ Borel measurable, say, with respect
  to the standard Hausdorff distance topology. The finiteness of the
  partition function ${\Bbb E}\sum_{\gamma \in \Gamma_D(\Pi_D)} \exp(-\Phi_D(\gamma))$
  above will be established in Section \ref{EoR}, Corollary \ref{JASNE}.


 \subsection{Dynamic representation}\label{DYKO}
  The availability of a dynamic description in terms of equilibrium time-space evolution
  of one-dimensional particle systems is a crucial feature of the Arak-Surgailis process and
  it underlies most of the techniques successfully used to derive explicit expressions
  for various characteristics of that model. It is therefore of particular importance
  to show that our three-dimensional model ${\cal A}_D$ admits an analogous dynamic
  representation in terms of evolution of two-dimensional edge systems tracing the
  polyhedral boundaries of ${\cal A}_D$ in three-dimensional time-space.

  Note that a planar section of a single polygonal face as defined above coincides,
  for $\mu$-almost all planes in ${\cal P},$ with a finite collection of disjoint
  colinear closed segments.
  Motivated by this observation we shall use the name of {\it multi-edge} for each
  such collection. For the purposes of this section we shall represent ${\Bbb R}^3$
  as three-dimensional time space ${\Bbb R} \times {\Bbb R}^2$ with two-dimensional
  spatial component; the time coordinate will be consistently denoted
  by $t.$ All the multi-edges considered below will arise as sections
  of polygonal faces in ${\Bbb R} \times {\Bbb R}^2$ with purely spatial planes $\varpi_t :=
  \{ t \} \times {\Bbb R}^2.$ In this setting, a non-empty multi-edge given as the
  intersection $\varepsilon := f \cap \varpi_t$ for some polygonal face
  $f,\; f \not\subseteq \varpi_t,$ can be assigned in a natural way its two-dimensional
  {\it velocity vector} lying in $\varpi_t$ and defined to be the spatial component
  (projection onto the spatial plane $\varpi_t$) of the
  unit normal to $\varepsilon$ lying in the plane $\pl[f]$ of $f$ and contained in
  the {\it future} half-space $\varpi_{t+} := [t,+\infty) \times {\Bbb R}^2,$
  divided by the length of the  spatial component of the unit normal to $\pl[f].$
  With $\varrho$ standing for the
  unique plane $\pl[f]$ containing the face $f$ we denote by $\vec{v}[\varrho]$
  the velocity vector constructed above (note that the construction does not
  depend on $t$). It is now easily verified that when observing the {\it time
  evolution} of the spatial plane $\varpi_t,$ we see the polygonal face $f$ mapped
  onto the time evolution of the multi-edge $\varpi_t \cap f$ {\it moving} in time
  with the velocity vector $\vec{v}[\pl(f)].$ Consequently, each finite family
  of disjoint and possibly nested bounded polyhedra is also mapped into the
  time evolution of the corresponding finite multi-edge systems. In course of
  the time evolution of multi-edges their components evolve as well, continuously
  updating their lengths according to the rule that no two edge segments can
  intersect each other except when sharing a vertex, whence the endpoints of a given edge
  segment are determined
  as the meeting points with other multi-edges. This rule may lead to extinction
  of certain edge segments as their length reaches $0$ but, on the other hand,
  certain segments may also get split into disjoint parts. A multi-edge dies
  if all its edges go extinct.

  As mentioned above, we shall consider our polyhedral random field
  defined in a certain open bounded and convex polyhedral domain
  $D \subseteq {\Bbb R} \times {\Bbb R}^2.$
  The presence of the boundary has its effect on the dynamics, which has to
  be taken into account in our construction. Namely, at each time
  $t \in {\Bbb R}$ we observe a collection of {\it boundary} multi-edges
  in $\varpi_t \cap \partial D$ and for each {\it internal} multi-edge
  in $D \cap \varpi_t$ the meeting point with the boundary marks the
  end of its appropriate segment. The boundary multi-edges can
  be assigned their velocity vectors on equal rights with the internal
  edges. It should be noted that a boundary multi-edge has always exactly
  one segment due to the convexity assumption imposed on $D.$

  Apart from the multi-edge extinction, the evolution of a multi-edge system,
  corresponding to a finite collection of disjoint and possibly nested
  polyhedral surfaces, comprises also multi-edge birth events. These can be
  divided into several groups
  \begin{description}
   \item{\bf (IT)} {\it Infinitesimal triangle birth} in a time-space point
         $x \in D \cap \varpi_t$ which does not lie on any internal or
         boundary multi-edge in $\varpi_t.$ At the point $x$ three
         pairwise non-parallel multi-edges $e_1, e_2, e_3$ are born,
         {\it moving} with three different velocities
         $\vec{v}_1, \vec{v}_2, \vec{v}_3.$
   \item{\bf (IA)} {\it Infinitesimal angle birth} in a time-space
         point $x$ lying on an internal or boundary edge $e_1$
         {\it moving} with a certain velocity $\vec{v}_1.$ At the point
         $x$ two multi-edges $e_2, e_3$ are born, {\it moving} with
         velocities $\vec{v}_2, \vec{v}_3$ and such that
         $\{ e_1, e_2, e_3 \}$ are pairwise non-parallel.
  \item{\bf (IE)} {\it Infinitesimal edge birth at a vertex} in a
        time-space point $x$ lying at the intersection of two different
        (either internal or boundary) edges $e_1, e_2 \subseteq
        (D \cup \partial D) \cap \varpi_t,$ {\it moving} with velocities
        $\vec{v}_1$ and $\vec{v}_2$ respectively. At the point $x$ a
        new multi-edge $e_3$ is born, {\it moving} with velocity $\vec{v}_3$
        and such that $\{ e_1, e_2, e_3 \}$ are pairwise non-parallel.
 \end{description}
  Observe that for certain choices of new-born multi-edges and their velocities we
  end up with unstable configurations where new edges die immediately upon
  their birth.
  This happens e.g. for an {\it infinitesimal triangle} with all its edges {\it moving}
  inward. We say that a multi-edge $e$ is {\it stably born} if it survives a positive
  amount of time rather than going extinct just upon its birth. The stability conditions
  are easily determined for all groups of birth events.
 \begin{description}
  \item{\bf (Stability for IT)} Denote by ${\bf n}_i$ the unit normal to $e_i$
        (in $\varpi_t$) pointing outward the infinitesimal triangle and let
        ${\bf v}[v_j,v_k]$ be the velocity vector (in $\varpi_t$) for the
        intersection point of $e_j$ and $e_k,$ easily checked to be
        $a v_j + b v_k,$ where
        $$ \left[ \begin{array}{c} a \\ b \end{array} \right]
        =  \left[ \begin{array}{cc} \langle v_j, v_j \rangle & \langle
                    v_j,v_k \rangle \\
                    \langle v_k,v_j \rangle & \langle v_k,v_k\rangle
                   \end{array} \right]^{-1}
           \left[ \begin{array}{c} |v_j|^2 \\ |v_k|^2 \end{array} \right]. $$ 
        The birth is stable iff $\langle \vec{v}_i, {\bf n}_i \rangle > 
        \langle {\bf v}[v_j,v_k], {\bf n}_i \rangle$
        for all $i \neq j \neq k \neq i.$ It can be easily shown that
        it is enough to verify this condition just for one fixed
        choice of $i,j,k.$
  \item {\bf (Stability for IA)} Let ${\bf n}_1$ be the unit normal to $e_1$
        in $\varpi_t$ pointing to the side of $e_1$ opposite to where
        the new angle is born and write ${\bf n}_2, {\bf n}_3$ for the unit
        normals to $e_2, e_3$ respectively, pointing outward the new-born
        convex angle. Once again, the birth is stable iff $\langle \vec{v}_i,
        {\bf n}_i \rangle  > \langle {\bf v}[v_j,v_k], {\bf n}_i \rangle$ for all
        $i \neq j \neq k \neq i.$
  \item {\bf (Stability for IE)} Let ${\bf n}_3$ be the unit normal to $e_3$ pointing inward
        the convex angle between $e_1$ and $e_2.$ Then the birth is stable iff
        $\langle \vec{v}_3, {\bf n}_3 \rangle > \langle {\bf v}[v_1,v_2], 
         {\bf n}_3 \rangle.$
  \end{description}
  In geometrical terms the stability condition for an {\it infinitesimal triangle}
  $e_1,e_2,e_3$ means that $e_i$ is able to {\it escape} in the direction pointed
  by ${\bf n}_i$ from the intersection point of $e_j$ and $e_k,$ whose
  velocity component in the direction ${\bf n}_i$ is smaller than
  $\langle \vec{v}_i, {\bf n}_i \rangle.$ Clearly, otherwise $e_i$ would be
  destroyed by other edges immediately upon its birth. The remaining stability
  conditions admit their geometrical interpretations along the same lines.

  In the sequel, we will refer to all stable {\it boundary} birth
  events (involving at least one boundary edge) of either {\bf (IA)}
  or {\bf (IE)} type as to {\it entry} events. Note that {\bf (IA)}
  or {\bf (IE)} events involving internal edges only are not
  considered to be entry events and neither are unstable birth
  events. For a given admissible polyhedral configuration $\gamma
  \in \Gamma_D$ we write $\entry(\gamma)$ for the collection
  of entry events it determines.

  Below, we construct a random dynamics on multi-edge
  systems conditionally on the collection ${\cal E}$ of all entry
  events assumed to be given. This stands in contrast to the
  two-dimensional Arak-Surgailis construction, where boundary birth events
  were also governed by explicit random dynamics, and this is due
  to the fact that we are unable to identify the explicit equilibrium
  distribution for entry events in the three-dimensional case.

  The random dynamics of internal birth events is given by
  the following rules, with $t$ standing for the time coordinate
  \begin{description}
   \item{\bf (Dynamics for IT)} The infinitesimal triangle birth
         time-space sites are chosen according to a homogeneous
         Poisson point process of intensity $\frac{\pi^4}{6}$ in $D.$
         The directions and velocities of the new-born edges are chosen
         according to the distribution of a typical stable vertex angle
         at a meeting point of three planes in $\Pi,$ where the
         adjective {\it stable} means here that out of the $8$ vertex
         angles at the triple intersection point we choose the {\it unique}
         one which gives rise to a stable infinitesimal triangle birth
         in the sense of {\bf (Stability for IT)},
   \item{\bf (Dynamics for IA)} The infinitesimal angle birth time-space
         sites on a face $f$ are chosen with intensity
         $\frac{\pi^3}{4} \sqrt{1+|v[\pl(f)]|^2} d\ell dt$
         with $d\ell$ standing for the edge length element on
         $f \cap \varpi_t$ (note that this intensity coincides
         with $\frac{\pi^3}{4}$ times the area element on $f$).
         The directions and velocities of the new-born edges are chosen
         according to the distribution of a typical stable vertex angle
         at a meeting point of $\pl(f)$  with two other planes of $\Pi,$
         where the adjective {\it stable} means here that out of the $8$
         vertex angles at the triple intersection point we choose the
         {\it unique} one which gives rise to a stable infinitesimal
         angle birth in the sense of {\bf (Stability for IA)},
   \item{\bf (Dynamics for IE)} An infinitesimal edge birth site
         at a vertex $x \in \varpi_t,$ at the intersection of two edges
         $e_1$ and $e_2$ lying respectively on faces $f_1$
         and $f_2$ and moving with respective velocities
         $\vec{v}_1$ and $\vec{v}_2,$ arises with intensity
         \begin{equation}\label{INTE1}
          \pi \sqrt{ 1+ |{\bf v}[v_1,v_2]|^2} d_2 x dt.
         \end{equation}
         Note that this intensity coincides with $\pi$ times the length element on
         the line traced by the time-space trajectory of the intersection
         point $e_1 \cap e_2.$ 
         The direction and velocity of the new-born edge are determined
         according to the distribution of a typical plane $\varpi$ of $\Pi$
         intersecting the line traced in time-space by $e_1 \cap e_2.$
         Unlike in previous cases, the stability of this birth event
         in the sense of {\bf (Stability for IE)}
         {\it is not guaranteed} here and in fact it is easily seen
         that it holds iff the normal to $\varpi$ pointing in time-space
         to the future
         is not contained in the solid angle between the faces $f_1$ and $f_2$
         [i.e. the solid angle between the half-planes of $\pl[f_1]$
         and $\pl[f_2]$ containing $e_1$ and $e_2$ respectively and
         meeting along $\pl[f_1] \cap \pl[f_2]$], denote this angle
         by $\angle(f_1,f_2)$ for use below. The non-stable birth
         events resulting from the above rule have no effect
         on the dynamics. Thus, alternatively we can produce the
         birth events with intensity
         \begin{equation}\label{INTE2}
          \pi \sqrt{ 1+ |{\bf v}[v_1,v_2]|^2}
           \frac{(2\pi-|\angle(f_1,f_2)|)}{2\pi} d_2 x dt
        \end{equation}
        and let the direction and velocity of the new-born edge be determined
        according to the distribution of a typical plane $\varpi$ of $\Pi$
        intersecting the line traced in time-space by $e_1 \cap e_2,$
        conditioned on yielding a stable infinitesimal edge birth event.
        In analogy with the interpretation of (\ref{INTE1}) above, the
        intensity (\ref{INTE2}) coincides with $\frac{1}{2} (2\pi - |\angle(f_1,f_2)|)$
        times the length element on the line $\pl[f_1] \cap \pl[f_2]$ traced
        by the time-space trajectory of the intersection point $e_1 \cap e_2.$
 \end{description}
  To make the above construction fully explicit, we note that the distribution
  of a typical vertex angle between three planes of $\Pi$ is $\propto
  |\langle {\bf n}_1,{\bf n}_2 \times {\bf n}_3 \rangle| d\sigma({\bf n}_1)
  d\sigma({\bf n}_2) d\sigma({\bf n}_3),$ where ${\bf n}_i \in {\Bbb S}_2,$
  $\sigma$ is the surface measure on ${\Bbb S}_2$ and the planes $\varpi_i,\;
  i=1,2,3$ creating the angle are respectively chosen orthogonal to
  ${\bf n}_i,\;i=1,2,3;$ indeed, this follows by a minor modification
  of Theorem 3 in Calka (2001), specialised for $d=3.$
  Recall that the scalar product $|\langle {\bf n}_1,{\bf n}_2 \times
  {\bf n}_3 \rangle|$ coincides with the volume of the parallelepiped
  spanned by the vectors ${\bf n}_1,{\bf n}_2, {\bf n}_3.$ Observe that
  the knowledge of the typical vertex angle distribution for three
  planes provides full knowledge of the corresponding laws with one
  or two planes fixed, as respectively required for {\bf (IA)} and
 {\bf (IE)}.

  The polyhedral process in $D$ resulting from the above dynamic construction
  in presence of a collection ${\cal E}$ of entry events will be denoted
  in the sequel by ${\cal A}_{D|{\cal E}}.$


\subsection{Equivalence of representations}\label{EoR}
 We now proceed to showing that both the Gibbsian representation
 (\ref{DEFI1}) and the dynamic construction of Subsection \ref{DYKO}
 yield, in a sense to be specified below, the same polyhedral
 field. To this end, for a given collection ${\cal E}$ of entry
 events in $D$ we put $\Gamma_{D|{\cal E}} := \{ \gamma \in
 \Gamma_D\;|\; \entry(\gamma) = {\cal E} \}.$
 Further, for a collection $\{ \varpi_i \}_{i=1}^n$
 of pairwise non-parallel planes in ${\Bbb R}^3$ we write
 $\Gamma_{D|{\cal E}}(\{\varpi\}_{i=1}^n)$ to denote the family of
 all admissible polyhedral configurations in $\Gamma_{D|{\cal E}}$
 with the additional property that
 $\{ \pl(f)\;|\; f \in F_D(\gamma) \setminus F[{\cal E}] \}$
 coincides with $\{ \varpi \}_{i=1}^n,$ where $F[{\cal E}]$ is the
 collection of planes arising in entry events from ${\cal E}.$
 Define $\kappa(D)$ to be the measure $\mu$ of the set of planes
 intersecting $D.$ Note that, $D$ being convex, we have
 \begin{equation}\label{KAPPA}
  \kappa(D) := \mu \{ \varpi,\; \varpi \cap D \neq \emptyset \}
  = \frac{1}{2} \sum_{e \in E(D)} |\angle(e)| \ell(e),
 \end{equation}
 see (4.2.30), (4.5.9) and (4.5.10) in Schneider (1993),
 whence $\kappa(D)$ is proportional to the generalised
 integral mean curvature functional of $D,$ see ibidem.

  The crucial observation in this Subsection is that, in view of Proposition \ref{WLASNOSCI},
  the form of the dynamic rules {\bf (Dynamics for IT,IA,IE)} implies that,
  given a collection ${\cal E}$ of entry events in $D,$ we have for all
  $\gamma \in \Gamma_{D|{\cal E}}$
  \begin{equation}\label{PRZEWAZNE}
   {\Bbb P}\left( {\cal A}_{D|{\cal E}} \in d\gamma \right)
  = \exp(\kappa(D)) {\bf 1}_{\{\gamma \in \Gamma_{D|{\cal E}}(\Pi_D) \}}
       \exp(-\Phi_D(\gamma)) d{\cal L}(\Pi_D)
 \end{equation}
 with ${\cal L}(X)$ standing for the law of a
 random object $X$. Indeed, for $\gamma \in \Gamma_{D|{\cal E}},$
 the probability differential $[d{\cal L}({\cal A}_{D|{\cal E}})](\gamma)$
 factorises into the product of the probability that all faces
 $f_1,\ldots,f_k$ of $F_D(\gamma) \setminus F[{\cal E}]$ were born in course of the
 evolution of
 the multi-edge system, which is seen to be $\prod_{i=1}^k [{\bf 1}_{\pl[f_i]
 \cap D \neq \emptyset} d\mu(\pl[f_i])] =
 \exp(\kappa(D)) [d{\cal L}(\Lambda_D)](\{ \pl[f_1],\ldots,\pl[f_k] \})$
 by comparing the dynamic rules {\bf (Dynamics for IT,IA,IE)} with Proposition
 \ref{WLASNOSCI}, times the probability that no other faces were born in the
 evolution giving rise to $\gamma,$ which is in its turn evaluated
 to $\exp(-\Phi_D(\gamma))$ in view of the form of the dynamic rules. This
 last claim is easily verified by noting that
 \begin{itemize}
  \item by {\bf (Dynamics for IT)}, $\exp\left(-\frac{\pi^4}{6} \vol(D)\right)$
        is the probability that no extra faces were born in {\bf (IT)} birth
        events,
  \item by {\bf (Dynamics for IA)}, $\exp\left(-\frac{\pi^3}{4} \sum_{f \in F_D(\gamma)} \Area(f) \right)$
        is the probability that no extra faces were born in {\bf (IA)} birth
        events,
  \item by {\bf (Dynamics for IE)}, $\exp\left(-\frac{1}{2} \sum_{e \in E_D(\gamma)}
        (2\pi - |\angle(e)|) \ell(e) \right)$ is the probability that no extra faces were born
        in {\bf (IE)} birth events, see the discussion following (\ref{INTE2}).
 \end{itemize}

 As an immediate consequence
 of (\ref{PRZEWAZNE}) we come to
 \begin{theorem}\label{ROWNOWAZNOSC}
  Given a collection ${\cal E}$ of entry events in $D,$ the
  random polyhedral field ${\cal A}_{D|{\cal E}}$ coincides
  in distribution with the Gibbs polyhedral field
  ${\cal A}^*_{D|{\cal E}}$ given by
  $$ {\Bbb P}\left({\cal A}^*_{D|{\cal E}} \in {\cal F}\right)
   = \frac{{\Bbb E}\sum_{\gamma \in \Gamma_{D|{\cal E}}(\Pi_D) \cap {\cal F}}
              \exp(-\Phi_D(\gamma))}
          {{\Bbb E}\sum_{\gamma \in \Gamma_{D|{\cal E}}(\Pi_D)}
              \exp(-\Phi_D(\gamma))} $$
  for each ${\cal F} \subseteq \Gamma_{D|{\cal E}}$ Borel measurable
  with respect to the usual Hausdorff topology. Moreover, for each
  ${\cal E}$ we have
  \begin{equation}\label{FUNKCJAPODZ}
   {\Bbb E}\sum_{\gamma \in \Gamma_{D|{\cal E}}(\Pi_D)}
           \exp(-\Phi_D(\gamma)) = \exp(-\kappa(D)).
  \end{equation}
 \end{theorem}
  It should be noted at this point that the particular form of
  the expression (\ref{FUNKCJAPODZ}) for the partition function
  might seem to stand in an unexpected contrast to the
  two-dimensional formula (4.6) in Arak \& Surgailis (1989).
  However, this difference is a matter of choice of the
  reference measure, which is in our case the normalised law
  of the Poisson plane process, while Arak \& Surgailis (1989)
  use $\sum_{n=0}^{\infty} \frac{1}{n!} \mu^{\otimes n}$
  instead. Moreover, we introduce the volume order term
  $\frac{\pi^4}{6} \vol(D)$ to our energy function,
  which has no equivalent in that paper.
  We find these choices preferable for the presentation
  of our setting, as leading to simpler formulae.
  However, should we use the reference measure and Hamiltonian
  analogous to those of that paper, our partition function in
  (\ref{FUNKCJAPODZ}) would evaluate to $\exp\left(\frac{\pi^4}{6} \vol(D)\right).$
  The absence of $\kappa(D)$ and a surface
  area order term in the exponent in this partition function expression
  is due to the fact that we {\it condition} on fixed collection of
  entry events rather than {\it randomising} it as done in
  Arak \& Surgailis (1989).

  As an obvious conclusion from Theorem \ref{ROWNOWAZNOSC} we get
  that the original polyhedral field ${\cal A}_D$ as defined by
  (\ref{DEFI1}) admits the dynamic representation with the collection
  of entry events distributed according to $\entry({\cal A}_D).$
  \begin{corollary}\label{JASNE}
   We have
   $$ {\cal L}({\cal A}_D) = \int {\cal L}({\cal A}_{D|{\cal E}})
                             d [{\cal L}(\entry({\cal A}_D))]({\cal E}) $$
   and
   $$ {\cal L}({\cal A}_{D|{\cal E}}) = {\cal L}({\cal A}_D|
      \entry({\cal A}_D)={\cal E}). $$
   Moreover,
   $$ {\Bbb E}\sum_{\gamma \in \Gamma_D(\Pi_D)} \exp(-\Phi_D(\gamma)) = \exp(-\kappa(D)). $$
  \end{corollary}
  The use of this corollary is limited by the fact that we do not
  know the distribution of $\entry({\cal A}_D).$

\section{Properties of the process}\label{WLAA}
 We argue first that the polyhedral fields ${\cal A}_D$ and
 ${\cal A}_{D|{\cal E}}$ exhibit a $3$-dimensional germ-Markov
 property in the sense specified in Corollary \ref{MARKOW3}
 below. For a smooth closed simple
 (non-self-intersecting) surface $\sigma$ in a bounded and
 convex polyhedral domain $D,$ by the trace of a polyhedral
 configuration $\gamma$ on $\sigma,$ denoted in the sequel
 by $\gamma \wedge \sigma,$ we mean the intersection $\gamma \cap \sigma$
 together with the directions of normals to the face planes at intersection.
 This concept can
 be formalised in various compatible ways, yet we keep the above
 informal definition in hope that it does not lead to any
 ambiguities while allowing us to avoid unnecessary technicalities.
 For convenience we assume that no face of $\gamma$ is tangent
 to $c,$ which can be ensured with probability $1$ in view
 of the smoothness of $\sigma.$
 In view of the Gibbsian representation (\ref{DEFI1}) and by
 Theorem \ref{ROWNOWAZNOSC} we easily conclude that
 \begin{corollary}\label{MARKOW3}
 For each $\sigma$ as above there exists
 a stochastic kernel ${\cal A}_{\innt \sigma}(\cdot|\vartheta)$
 such that, with $\vartheta$ standing for a trace on $\sigma,$
 \begin{equation}\label{WAARUNK}
  {\cal L}_{\innt \sigma}\left({\cal A}_D | {\cal A}_D \wedge
    \sigma = \vartheta \right) =
  {\cal L}_{\innt \sigma}\left({\cal A}_{D|{\cal E}} |
         {\cal A}_{D|{\cal E}} \wedge
    \sigma = \vartheta \right) = {\cal A}_{\innt \sigma}(\cdot|\vartheta)
 \end{equation}
 for all bounded open and convex polyhedral domains $D \supseteq \overline{\innt \sigma}$
 and for each collection ${\cal E}$ of entry events in $D,$ where
 ${\cal L}_{\innt \sigma}(X)$ denotes the law of a
 random field $X$ restricted to $\innt \sigma$
 (the interior of $\sigma$).
 \end{corollary}
 Indeed, it is easily seen that we have for
 measurable $G \subseteq \Gamma_{\innt \sigma||\vartheta}$
 \begin{equation}\label{OPISJADRA}
  {\cal A}_{\innt \sigma}(G |\vartheta) =
  \frac{{\Bbb E} \sum_{\gamma \in \Gamma_{\innt \sigma || \vartheta}(\Pi_{\innt \sigma}) \cap G}
        \exp(-\Phi_{\innt \sigma}(\gamma))}
       {{\Bbb E} \sum_{\gamma
        \in \Gamma_{\innt \sigma || \vartheta}(\Pi_{\innt \sigma})} \exp(-\Phi_{\innt \sigma}(\gamma))},
 \end{equation}
 where $\Gamma_{\innt \sigma || \vartheta} := \{ \gamma \in
 \Gamma_{\innt \sigma}\;|\; \gamma \wedge \sigma = \vartheta \}$
 and, for a collection $\{ \varpi_i \}_{i=1}^n$ of planes in ${\Bbb R}^3$
 hitting $\innt \sigma,$ $\Gamma_{\innt \sigma || \vartheta}(\{ \varpi_i \}_{i=1}^n)$
 denotes the family of all polyhedral
 configurations $\gamma \in \Gamma_{\innt \sigma || \vartheta}$
 for which the set $\{ \pl[f]\;|\; f \in F_D(\gamma),\; f \cap \sigma =
 \emptyset \}$ coincides with $\{ \varpi_i \}_{i=1}^{n}.$

 To proceed, consider the family $\Gamma_{{\Bbb R}^3}$ of whole-space
 admissible polyhedral configurations, determined by {\bf (P1), (P2), (P4)}
 and {\bf (P6)} ({\bf (P3)} and {\bf (P5)} are meaningless in this context)
 and by the requirement of local finiteness (any bounded set is hit by at
 most a finite number of faces). It is natural to define the family
 ${\cal G}({\cal A})$ of infinite volume Gibbs measures (thermodynamic limits)
 for ${\cal A}$ as the collection of all probability measures on
 $\Gamma_{{\Bbb R}^3}$ with the accordingly distributed random element
 ${\cal A}$ satisfying
 \begin{equation}\label{GRANTERM}
  {\cal L}_{\innt \sigma}\left({\cal A} | {\cal A} \wedge \sigma = \vartheta \right) =
  {\cal A}_{\innt \sigma}(\cdot|\vartheta)
 \end{equation}
 for $\sigma$ ranging through the collection of all bounded smooth
 simple closed surfaces in ${\Bbb R}^3.$
 In addition, we shall consider the family ${\cal G}_{\tau}({\cal A})$
 of isometry invariant measures in ${\cal G}({\cal A}).$ We believe
 that, in analogy with the results in Section 3 in Schreiber (2005),
 it should be possible to show that ${\cal G}_{\tau}({\cal A}) \neq \emptyset$
 by using an appropriate relative compactness argument. Moreover, we
 conjecture that the uniqueness of the isometry invariant thermodynamic
 limit ${\cal A}_{{\Bbb R}^3},$ as well as its coincidence with the polyhedral
 process traced by infinite-volume equlibrium evolution of the multi-edge
 system as discussed above, could possibly be established following the
 lines of Schreiber (2004a), where this is done for two-dimensional
 polygonal fields admitting dynamic representation.
 On the other hand, in analogy with the two-dimensional setting, we do
 not expect that ${\cal A}_{{\Bbb R}^3}$ be the unique element of
 ${\cal G}({\cal A}),$ see the discussion closing Section 3 in
 Schreiber (2005). We are working on this conjecture at present,
 yet we are unable to provide their formal proofs
 at the current stage of our research. However, should these conjectures
 hold true as stated, initiating the multi-edge system dynamics in a
 domain $D$ with the collection of entry events
 $\entry({\cal A}_{{\Bbb R}^3} \cap D),$
 for ${\cal A}_{{\Bbb R}^3}$ denoting the unique thermodynamic limit,
 would result in a consistent family of polyhedral fields, as constructed
 by Arak \& Surgailis (1989) for the polygonal setting. A further
 essential task would be to provide a feasible description of the
 entry process $\entry({\cal A}_{{\Bbb R}^3} \cap D),$ which we
 anticipate to be of a rather complicated nature.


\section{Birth site birth and death dynamics for simulating
         polyhedral fields}\label{BSDYN}
 The purpose of the current section is to construct, much along the
 lines of Schreiber (2005), Section 2.1, a random dynamics on the space
 $\Gamma_D$ of admissible polyhedral configurations which
 leaves invariant the law of ${\cal A}_{D|{\cal E}},$ where
 the collection ${\cal E}$ of entry events is to remain fixed
 throughout the section. This will allow us later to provide
 modifications of this dynamics suitable for simulation of
 Gibbsian modifications of ${\cal A}_{D|{\cal E}}.$
 The purpose of developing this algorithm is its envisioned
 application, as a component of suitable simulated annealing
 techniques, to volumetric image segmentation along the lines
 of our previous papers Kluszczy\'nski et al. (2004,2005)
 where we considered the corresponding two-dimensional problem.

 In the sequel, particular care is needed to distinguish between the notion
 of time considered in the dynamic representation of the polygonal field
 ${\cal A}_{D|{\cal E}}$ given in Subsection \ref{DYKO} above, and the notion of
 time to be introduced for the random dynamics on $\Gamma_{D|{\cal E}}$
 constructed below. To make this distinction clear we shall refer to the
 former as to the {\it representation time} (r-time for short) and shall
 keep for it the notation $t,$ while the latter will be called the
 {\it simulation time} (s-time for short) and will be consequently
 denoted by $s$ in the sequel.

 It is convenient for our exposition below to perceive each individual
 infinitesimal triangle birth site ({\bf (IT)}-birth site) in the dynamic
 representation, see {\bf (IT)}, as coming with an associated random number
 generator, represented for instance as an infinite sequence of i.i.d. random
 variables uniformly distributed on $[0,1]$ and used to determine
 the subsequent moments and angles$\slash$velocities for critical
 events {\bf (IA)} and {\bf (IE)} involving multi-edges resulting
 from the considered {\bf (IT)} birth event. In other words, each
 {\bf (IT)}-birth site is assumed to carry a {\it package} enclosing all
 randomness the resulting multi-edges may possibly encounter during their evolution,
 and the above is just one technical possibility of how this can be
 achieved. We shall use the name of a {\it birth package} for an
 infinitesimal triangle birth site with such a {\it random number
 generator} attached. In these terms, it is now easily seen that
 the polyhedral configuration obtained in course of the dynamic
 construction depends {\it deterministically} on the underlying
 collection of birth packages.

 Consider a polyhedral configuration $\gamma \in \Gamma_{D|{\cal E}}$ and
 a new infinitesimal triangle birth site $x_0 \in D$ not yet present
 in $\gamma,$  extended to a birth package in the standard way as
 discussed above. Adding this birth package to the collection of birth
 packages determining $\gamma$ and keeping the evolution
 rules of the dynamic representation {\bf (Dynamics for IT,IA,IE)}
 results in a new configuration to be denoted by $\gamma \oplus x_0.$
 Likewise, removing an {\bf (IT)}-birth site $x_1$ from a configuration
 $\gamma$ in which it was present yields a new polyhedral configuration
 $\gamma \ominus x_1.$

 Taking into account that the collection of the ${\bf (IT)}$-birth sites
 for ${\cal A}_{D|{\cal E}}$ is chosen according to the Poisson point process
 with intensity $\frac{\pi^4}{6}$ as specified in {\bf (Dynamics
 for IT)}, we easily see that the law of ${\cal A}_{D|{\cal E}}$
 is invariant with respect to the following pure-jump Markovian
 {\bf (IT)}-birth site birth and death dynamics on $\Gamma_{D|{\cal E}},$
 further denoted by {\bf (BS)}, with $\gamma_s$ standing for the
 state at time $s$ and with $\gamma_0 \in \Gamma_{D|{\cal E}}$
 \begin{description}
  \item{\bf (BS:birth)} With intensity $\frac{\pi^4}{6}ds$ set
         $\gamma_{s+ds} := \gamma_s \oplus x,$
  \item{\bf (BS:death)} For each {\bf (IT)}-birth site in $\gamma_s$ with
         intensity $1 \cdot ds$ set
         $\gamma_{s+ds} := \gamma_s \ominus x,$
  \item{}  otherwise keep $\gamma_{s+ds} = \gamma_s.$
  \end{description}
  In fact, more can be stated, see also Proposition 1 in Schreiber (2005)
 \begin{theorem}\label{SIM1}
  The distribution of the polygonal field ${\cal A}_{D|{\cal E}}$ is
  the unique invariant law of the dynamics given by {\bf (BS:birth)},
  and {\bf (BS:death)}. The resulting stationary
  process is reversible. Moreover, for any initial distribution of
  $\gamma_0$ concentrated on $\Gamma_{D|{\cal E}}$
  the laws of the random polygonal fields $\gamma_s$
  converge in total variation to the law of ${\cal A}_{D|{\cal E}}$
  as $s \to \infty.$
 \end{theorem}
 While the invariance was discussed above and the reversibility is
 clear, the uniqueness and convergence statements in the above theorem
 require a short justification. They both follow by the observation
 that, in finite volume, regardless of the initial state, the process
 $\gamma_s$ spends a non-null fraction of time in the state where
 no polyhedral faces other than those arising in ${\cal E}$ are present.
 Indeed, this observation allows us to conclude the required uniqueness
 and convergence by a standard coupling argument.

 \paragraph{Dynamics for Gibbsian modifications}
  Assume that a Hamiltonian (energy function) ${\cal H}$ is defined
  on the space $\Gamma_D$ of admissible polyhedral
  configurations and it satisfies
  \begin{equation}\label{HAMILTONIAN}
   {\cal H}(\gamma) \geq - A \vol(D) - B
  \end{equation}
  for some positive constants $A,B.$ Then it is clear
  that the partition function
  $$ Z_D[{\cal H}] := {\Bbb E}\exp(-{\cal H}({\cal A}_{D|{\cal E}})) $$
  is finite. Consequently,
  the corresponding Gibbsian modification ${\cal A}_{D|{\cal E}}^{\cal H}$
  can be considered with
  \begin{equation}\label{GibbsModified}
   \frac{d {\cal L}({\cal A}_{D|{\cal E}}^{\cal H})}
        {d {\cal L}({\cal A}_{D|{\cal E}})}[\gamma]
    = \frac{\exp(-{\cal H}(\gamma))}{Z_D[{\cal H}]},\; \gamma \in \Gamma_{D|{\cal E}}.
  \end{equation}
  Consider the following modification of the basic ${\bf (BS)}$
  dynamics:
  \begin{description}
   \item{\bf (BS[{\cal H}]:birth)}
         With intensity $\frac{\pi^4}{6}ds$ propose the update
         $\delta := \gamma_s \oplus x.$ Then,
         with probability $\min(1,\exp({\cal H}(\gamma_s)
         -{\cal H}(\delta)))$ accept this update, putting
          $\gamma_{s+ds} := \delta,$ otherwise keep $\gamma_{s+ds}
         := \gamma_s,$
  \item{\bf (BS[{\cal H}]:death)} For each {\bf (IT)}-birth site
         in $\gamma_s$ with intensity $1 \cdot ds$ set
         $\delta := \gamma_s \ominus x.$
         Then, with probability $\min(1,\exp({\cal H}(\gamma_s)
         -{\cal H}(\delta)))$ put $\gamma_{s+ds}
         := \delta,$ otherwise keep $\gamma_{s+ds}
         := \gamma_s,$
 \end{description}
  In other words, the original dynamics {\bf (BS)} is used in the
  standard way to propose a new configuration $\delta,$ which
  is then accepted with probability  $\min(1,\exp({\cal H}(\gamma_s)
  -{\cal H}(\delta)))$ and rejected otherwise.
  As a direct consequence of Theorem \ref{SIM1} and by a standard
  check of the detailed balance condition we get
 \begin{theorem}\label{SIM2}
   The distribution of the polyhedral field ${\cal A}^{\cal H}_{D|{\cal E}}$ is
   the unique invariant law of the dynamics given by
   {\bf (BS[{\cal H}]:birth)} and
   {\bf (BS[{\cal H}]:death)}. The resulting
   stationary process is reversible. Moreover, for any initial distribution of
   $\gamma_0$ concentrated on $\Gamma_{D|{\cal E}}$
   the laws of the random polyhedral fields $\gamma_s$
   converge in total variation to the law of ${\cal A}^{\cal H}_{D|{\cal E}}$ as
   $s \to \infty.$
 \end{theorem}
 \paragraph{Acknowledgements}
   The author gratefully acknowledges the support of the {\it Foundation for
   Polish Science} (FNP), from the Polish Minister of
   Scientific Research and Information Technology grant 1 P03A 018 28 (2005-2007)
   and from
   the EC 6th Framework Programme Priority 2 Information Society Technology
   Network of Excellence MUSCLE (Multimedia Understanding through Semantics,
   Computation and Learning; FP6-507752). He also wishes to express his gratitude
   for hospitality of Marie-Colette van Lieshout at the {\it Centrum voor Wiskunde
   en Informatica} (CWI), Amsterdam, The Netherlands, where a part of this research
   was completed. Thanks are also due to Pierre Calka for his valuable comments
   and suggestions.

 \paragraph{References}
 \begin{description}
  \item{\sc Arak, T.} (1982) On Markovian random fields with finite number of values,
       {\it 4th USSR-Japan symposium on probability theory and mathematical statistics,
           Abstracts of Communications}, Tbilisi.
  \item{\sc Arak, T., Surgailis, D.} (1989) Markov Fields with Polygonal Realisations,
       {\it Probab. Th. Rel. Fields} {\bf 80}, 543-579.
  \item{\sc Arak, T., Surgailis, D.} (1991) Consistent polygonal fields,
       {\it Probab. Th. Rel. Fields} {\bf 89}, 319-346.
  \item{\sc Arak, T., Clifford, P., Surgailis, D.} (1993) Point-based polygonal models
        for random graphs, {\it Adv. Appl. Probab.} {\bf 25}, 348-372.
  \item{\sc Calka, P.} (2001) Mosa\"iques poisoniennes de l'espace euclidien. Une
        extension d'un r\'esultat de R.E. Miles, {\it C. R. Acad. Sci. Paris, S\'er. I Math.},
        {\bf 332}(6) 557-562.
  \item{\sc Clifford, P., Middleton, R.D.} (1989) Reconstruction of polygonal images,
        {\it Journal of Applied Statistics}, {\bf 16}, 409-422.
  \item{\sc Clifford, P., Nicholls, G.} (1994) A Metropolis sampler
        for polygonal image reconstruction, {\it available at}:\\
        {\tt http://www.stats.ox.ac.uk/~clifford/papers/met\_poly.html},
 \item {\sc Kluszcz\'ynski, R., van Lieshout, M.N.M. and Schreiber, T.} (2004)
        Image segmentation by polygonal Markov fields, {\it submitted},
        Electronic version available as a CWI Research Report PNA-R0409 at:
        {\tt http://www.cwi.nl/publications}.
  \item{\sc Kluszczy\'nski, R., van Lieshout, M.N.M. and Schreiber, T.} (2005)
        An algorithm for binary image segmentation using polygonal Markov fields,
        accepted for the conference ICIAP (2005) [International Conference
        on Image Analysis and Processing, Cagliari, Italy], to appear
        in the proceedings (LNCS).
  \item{\sc Matheron, G.} (1975) {\it Random sets and integral geometry}, Wiley \& Sons,
        New York.
  \item{\sc Miles, R.E.} (1971) Poisson flats in Euclidean spaces, Part II: Homogeneous
         Poisson flats and the complementary theorem. {\it Adv. Appl. Probab.} {\bf 3},
         1-43.
 \item{\sc Schneider, R.} (1993) {\it Convex bodies: The Brunn-Minkowski Theory}(Encyclopaedia
         Math. Appl. {\bf 44}), Cambridge University Press.
  \item{\sc Schreiber, T.} (2004a) Mixing properties for polygonal Markov fields in the plane,
        {\it submitted}, available at: {\tt http://www.mat.uni.torun.pl/preprints},
  \item {\sc Schreiber, T.} (2005).
         Random dynamics and thermodynamic limits for polygonal Markov 
         fields in the plane.
         {\it Adv. in Appl. Probab.} {\bf 37}, 884-907.
 \end{description}

\end{document}